\documentclass[11pt]{amsart}
\usepackage{amsthm}
\usepackage{amsmath}
\usepackage{amsfonts}
\usepackage{amssymb}
\usepackage{color}
\newcommand{\ind}{{\rm ind}}

\newcommand{\beq}{\begin{equation}}
\newcommand{\beqn}{\begin{equation*}}

\newcommand{\eeq}{\end{equation}}
\newcommand{\eeqn}{\end{equation*}}

\newcommand{\R}{\mathbb{R}}

\newcommand{\LS}{\Lambda S^n}
\newcommand{\LSzero}{\Lambda^0 S^n}
\newcommand{\n}{\mathbb{N}}

\newcommand{\q}{\mathbb{Q}}

\newcommand{\sn}{S^n}
\newcommand{\sone}{S^1}
\newcommand{\z}{\mathbb{Z}}
\newtheorem{theorem}{Theorem}

\newtheorem{proposition}{Proposition}

\title[bumpy metrics and minimal
index growth]{Bumpy metrics on spheres and
minimal index growth}
\author{Hans-Bert Rademacher}
\address{Mathematisches Institut, Universit{\"a}t Leipzig,
04081 Leipzig, Germany}
\email{rademacher@math.uni-leipzig.de}
\date{2016-08-04, revised: 2016-09-02}
\dedicatory{Dedicated to Paul Rabinowitz with best wishes
}
\begin{document}
\begin{abstract}
The existence of
two geometrically distinct closed geo\-desics on an
$n$-dimensional sphere $S^n$ with a non-reversible 
and bumpy Finsler metric
was shown
independently 
by Duan \& Long~\cite{DL} and the author~\cite{Ra10}.
We simplify the proof of this statement by 
the following observation: If 
for some $N \in \n$
all closed geodesics 
of index $\le N$
of a non-reversible and bumpy
Finsler metric on $S^n$
are geometrically equivalent 
to the closed geodesic $c$ then
there is a covering $c^r$
of minimal index growth, i.e.
$$\ind(c^{rm})=m \,\ind(c^r)-(m-1)(n-1)$$
for all $m \ge 1$ 
with $\ind\left(c^{rm}\right)\le N.$
But this leads to a contradiction 
for $N =\infty$
as pointed out
by Goresky \& Hingston~\cite{GH}.
We also discuss perturbations of Katok metrics
on spheres of even dimension
carrying only finitely many closed geodesics. 
For
arbitrarily large $L>0$ we obtain on $S^2$ a metric 
of positive flag curvature
carrying only two closed geodesics
of length $<L$ which do not intersect.
\end{abstract}
\maketitle
\section{Introduction}
In this note we consider existence results for closed
geodesics on spheres $\sn$ 
of dimension $n$
endowed with a 
Finsler metric $f.$
The Finsler metric $f$ is called {\em reversible,} if 
$f(-X)=f(X)$ holds for all tangent
vectors $X.$ Otherwise we call the metric
{\em non-reversible.} Note that two closed geodesics $c_1,c_2: \sone \longrightarrow M$
of a non-reversible Finsler metric are 
{\em geometrically equivalent}
if their images $c_1(\sone)=c_2(\sone)$ 
and their orientations coincide. In the reversible case $c_1,c_2$
are geometrically equivalent if 
and only if $c_1(\sone)=c_2(\sone).$
Closed geodesics which are not geometrically equivalent are
called {\em geometrically distinct.}
We call a closed geodesic {\em prime} if it is not the covering of a shorter
closed curve.
For $m\ge 1$ and a prime closed geodesic $c$ of a non-reversible
Finsler metric the coverings $c^m$ defined
by $c^m(t)=c(mt)$ are geometrically
equivalent closed geodesics. 

A closed geodesic $c: \sone \longrightarrow \sn$
is {\em non-degenerate}
if its nullity vanishes.
The nullity is the 
dimension of the kernel of the index form minus one, it
equals the dimension of periodic Jacobi fields orthogonal 
to the velocity field $c'$ of the closed geodesic.

A Finsler metric $f$ on a compact manifold $M$ is 
called {\em bumpy} if all 
closed geodesics are non-degenerate.
In this case the energy functional 
$E:\Lambda \sn \longrightarrow \R$
defined on the free loop space $\LS$
carrying a canonical $\sone$-action 
can be viewed as a Morse function whose critical set
decomposes into a 
disjoint union of non-degenerate $S^1$-orbits 
of closed geodesics.
Following ideas by Birkhoff~\cite[Sec.17]{B1917},
~\cite[p.135-139]{B1927} the existence of a single closed
geodesic on a compact and simply-connected manifold with
a Finsler metric was shown by Lusternik \& Fet~\cite{LF}.
It turns out that existence results for several closed
geodesics strongly depend on the reversibility of the metric.
The sequence $\left(\ind\left(c^m\right)\right)_{m\ge 1}$ of Morse indices of
the coverings $c^m$ of a closed geodesic $c$
plays an important
role in existence proofs for several closed geodesics.
It is a result by Fet~\cite{Fe} that a bumpy and reversible Finsler metric
on a compact and simply-connected manifold carries at least two 
geometrically distinct closed geodesics.
An analogous statement holds also for non-reversible Finsler metrics on spheres
as shown independently by Duan-Long in~\cite{DL} and
the author in~\cite{Ra10}:
\begin{theorem}
\label{thm:eins}
{\rm \cite{DL}, \cite{Ra10}}
A bumpy and non-reversible Finsler metric on a sphere $S^n$ of dimension 
$n \ge 3$ carries
two geometrically distinct closed geodesics.
\end{theorem}
In~\cite{DL} the classification
of the 
symplectic normal forms of the
linearized Poincar\'e mapping of a closed geodesic
and a case distinction
is  used. In~\cite{Ra10} the {\em common index jump
theorem} due to~\cite[Thm.4.3]{LZ} is the main
ingredient.
The statement of Theorem~\ref{thm:eins} as well as the
statement of the following 
Proposition~\ref{pro:eins} also hold
for compact and simply-connected manifolds $M$
which are rationally homotopy equivalent to an
$n$-sphere $\sn.$
In this note we show that
one can use the following result 
about the index growth, i.e. the growth of the 
sequence $\left(\ind\left(c^m\right)\right)_{m\ge 1}$
which also gives a short proof of Theorem~\ref{thm:eins}:
\begin{proposition}
\label{pro:eins}
Let $f$ be a bumpy and non-reversible 
Finsler metric on $\sn$
and let $N_1 \ge 5n$ be a number such that the following
assumption holds:
There is a prime closed geodesic $c$ such that all closed
geodesics $d$ whose index satisfies $\ind(d)\le N_1+2$
are geometrically equivalent to $c.$
Let $r=n$ if $n$ is even and $r=(n+1)/2$
if $n$ is odd. Let 
$m_1:=\max\left\{m\in \n\,;\,\ind(c^{rm})\le N_1 \right\}\,.$
Then $m_1\ge 2$ and 
the closed geodesic $c^r$ is {\em of minimal index growth
up to level $m_1,$}
i.e.
\begin{equation}
 \label{eq:minimalequality}
 \ind(c^{rm})= m\,\ind(c^r)-(m-1)(n-1)
 \end{equation}
holds for all $m \le m_1.$
The closed geodesic $c$ is of
{\em elliptic-parabolic type,} i.e. its linearized Poincar\'e
mapping decomposes into rotations.
\end{proposition}
We prove this Proposition in the next section.
In the case $n=2, N_1=\infty$ 
it was shown by Ziller in~\cite[p.149]{Zi} that 
Equation~\eqref{eq:minimalequality} holds
with $r=2$ and leads to a contradiction.

The closed geodesic $c^r$ satisfying
Equality~\eqref{eq:minimalequality} is called 
{\em of minimal index growth up to level 
$m_1$}. This is motivated by the following 
Inequality~\eqref{eq:minimalinequality}.
On the other hand the following
Inequality~\eqref{eq:minimalplus} shows that
Equation~\eqref{eq:minimalequality}
can hold only for finitely many $m$
which shows that Proposition~\ref{pro:eins} implies
Theorem~\ref{thm:eins} for $N_1=\infty.$
\begin{proposition}
\label{pro:zwei}
{\rm (cf.~Goresky \& Hingston~\cite[Prop.6.1]{GH})}
Let $c$ be a closed geodesic of a bumpy Finsler metric.
Then for all $m\ge 1:$ 
\begin{equation}
 \label{eq:minimalinequality}
 \ind(c^{rm})\ge m\,\ind(c^r)-(m-1)(n-1)\,.
 \end{equation}
And there is an $N \in \n$ such that for all
$m > N:$
\begin{equation}
 \label{eq:minimalplus}
 \ind(c^{rm})\ge m\,\ind(c^r)-(m-1)(n-1)+1\,.
 \end{equation}
\end{proposition}
This Proposition is a particular case of~\cite[Prop.6.1]{GH}.
It is a direct consequence of Bott's formula for
$\ind(c^m),$ cf.~\cite[Thm.A]{Bo}.
The statement of the Proposition can also be shown without the
assumption that the metric is bumpy. 

It is a remarkable result by Bangert \& Long~\cite {BL}
that for {\em any} non-reversible Finsler metric on $S^2$ there are
two closed geodesics. There are a number of results for 
metrics in higher dimensions with and without non-degeneracy
assumptions or with curvature assumptions,
cf. for example
\cite{DL09},\cite{DL10a}, \cite{DL10b},
\cite{DLW15}, \cite{HR} and \cite{Ra07}.

There are non-reversible bumpy Finsler metrics 
of constant flag curvature
depending on an irrational parameter
on
$S^{2n-1}$ resp. 
$S^{2n}$ with exactly $2n$  geometrically distinct
closed geodesics.
On the $2$-sphere for any $N_1 \in \n$ there is a 
Katok metric
satisfying the assumptions of 
Proposition~\ref{pro:eins}.
These metrics occured in a paper by Katok~\cite{Ka},
their geometry is investigated by Ziller in~\cite{Zi}. 
By construction these metrics are invariant under 
a rotation,
the closed geodesics 
$c_1,c_1^{-1}, \ldots,c_{n},c_n^{-1}$
occur in pairs differing by orientation and length.
They 
are fixed point sets of 
reflections at a two-dimensional plane.
In Section~\ref{sec:example}
we show that for a given
arbitrarily large $L$
one can perturb 
these metrics on $S^{2n}$ to obtain a
non-reversible Finsler metric
with positive flag curvature
carrying exactly $2n$
geometrically distinct closed geodesics
$d_1,d_2,\ldots,d_{2n}$
with length $<L,$
which pairwise do not intersect.
Hence for these closed geodesics 
the closed curves $d_j^{-1}, j=1,2,\ldots,2n$ 
with opposite orientation
are not geodesics.
These metrics  are also invariant under a rotation
and are obtained by breaking the symmetry with respect to 
the reflections at two-dimensional planes.

In dimension $2$ this produces for $k \ge 4$
in any $C^k$-neighborhood of the Katok metric
examples of non-reversible Finsler metrics
on $S^2$ of positive
flag curvature with two simple closed geodesics, which do
not intersect.
\section*{Acknowledgement}
I am grateful to Nancy Hingston for helpful discussions
about the topic of the paper and her comments on an earlier version.
And I want to thank Philip Kupper for his careful reading and
the referee for his suggestions.
\section{Proof of Proposition~\ref{pro:eins}}
We denote by $f: T\sn \longrightarrow \R$ 
the Finsler metric defined on the tangent bundle
$T\sn$ of the $n$-sphere. 
We consider the free loop space $\Lambda S^n$ of $H^1$-maps 
$\sigma:S^1=[0,1]/\{0,1\}\longrightarrow S^n$
on the sphere $S^n.$ 

Closed geodesics are the critical points of the energy functional
$$E: \Lambda S^n \longrightarrow \mathbb{R}\,;\,
E(\sigma)=\frac{1}{2}\int_{0}^1 f^2(\sigma'(t))\,dt.$$
There is a canonical $S^1$-action 
$
(z,\sigma)\in S^1\times \Lambda S^n
\mapsto z.c \in \LS$
defined by $z.c(t)=c(t+z)$
on the free loop space
leaving the energy functional invariant.
The index form $\mathcal{H}_c$ of the closed geodesic 
can be identified with the hessian 
$d^2 E(c)$ of the energy functional
at $c$ by the second variational formula.
If the metric is bumpy one can view the energy functional
as a Morse function on the quotient space $\LS/\sone,$
for details see for example~\cite[Sec.2]{Ra10}.

We assume that there is a prime closed geodesic $c$ such that any other
closed geodesic $d$ with $\ind(d)\le N_1+2$ satisfies $d=z.c^m$
for some $z\in S^1$ and $m\ge 1.$ 

Let $\LSzero=\{\sigma \in \LS\,;\,E(\sigma)=0\}$
be the set of point curves, which one can identify with
the sphere $S^n.$ We denote by 
$b_k:=b_k\left(\LS/\sone,\LSzero/\sone;\q\right)$
the Betti numbers of the quotient space pair 
$\left(\LS/\sone, \LSzero/\sone\right)$ with rational
coefficients. 

Let $\Lambda(c):=\{\sigma \in \Lambda S^n\,;\, E(\sigma)<E(c).$
Then the Betti numbers of the 
local homology produced by
the covering $c^m$ are given by:
\begin{eqnarray*}
b_j(c^m)&=& \dim H_j\left(\left(\Lambda(c^m)\cup \sone.c^m\right)/\sone,
\Lambda(c^m)/\sone;\q\right)\\
&=&
\left\{
\begin{array}{ccl}
1&;& j=\ind(c^m) \mbox{ \em and }\\
&& m \equiv 1 \mbox{ \em or }
\ind(c^2)\equiv \ind(c) \pmod{2}\\
0 &;& otherwise
\end{array}\,.
\right.
\end{eqnarray*}
For $k \le N_1+2:$
\begin{equation}
 \label{eq:wk}
 w_k=\# \{m \in \n\,;\,
  \ind(c^m)=k, m \equiv 1 \mbox{ \em or }
  \ind(c^2)\equiv \ind(c) \pmod{2}\}
\end{equation}
gives the number of critical points 
of index $k$ of the Morse function
$E: \LS/\sone \longrightarrow \R$ 
producing non-trivial local homology.

Bott's formula for the 
indices $\ind(c^m)$ of iterates
$c^m$ implies the following statements
about the parity of $\ind(c^m), m\ge 1:$
\begin{equation}
\label{eq:parity}
\ind(c^{2m+1})\equiv \ind(c) \,;\, 
\ind(c^{2m})\equiv \ind(c^2)\pmod{2}\,,
\end{equation}
cf.~\cite[Sec.1]{Ra89}.
This observation implies that
$w_k=0$ for all numbers $k \equiv \ind(c)+1 \pmod{2}$ 
for all $k\le N_1+2.$
Hence 
the equality case holds in the Morse inequalities:
\begin{equation}
\label{eq:morse}
w_k=b_k, k\le N_1\,,
\end{equation}
cf.~\cite[(2.3)]{Ra89}.
Since Bott's formula also implies
$\ind(c^m)\ge \ind(c)\,;\, m\ge 1$
we conclude 
\begin{equation}
\label{eq:ceins}
\ind(c)=n-1\,.
\end{equation}
Here we use that $b_0=b_1=\ldots=b_{n-2}=0, b_{n-1}=1,$
cf. \cite[Thm.2.4]{Ra89} resp.
\cite[p.104]{Hi84}.
Equation~\eqref{eq:ceins} implies that the sequence
$\left(\ind(c^m)\right)_{m\ge 1}$ is monotone increasing,
i.e.
\begin{equation}
\label{eq:monotone}
\ind(c^{m+1})\ge \ind(c^m)\,,\, m\ge 1\,,
\end{equation}
cf. the {\em successive index estimates}
by Long \& Zhu, cf.~\cite{LZ}, 
\cite[Sec. 10.2]{Lo}.
If $\ind(c^2)\equiv \ind(c)+1 \pmod{2}$ it
follows from Inequality~\eqref{eq:monotone}
and Equality~\eqref{eq:parity} that
\begin{equation}
\label{eq:monotone1}
\ind(c^{m+1})\ge \ind(c^m)+1
\end{equation}
for all $m\ge 1.$

Now we discuss four cases depending on the 
parity of the dimension $n$
and the parity of $\ind(c^2)-\ind(c).$

\vspace{0.7cm}

{\sc Case 1:} Assume $n$ even, 
then
 \begin{equation}
  b_k=
  \left\{
  \begin{array}{ccl}
  2 &;& k=(2j+1)(n-1), j\ge 1\\
  1&;& k\ge n-1, k \mbox{ odd }, k \not=(2j+1)(n-1), j\ge 1\\
  0&;& \mbox{ otherwise }
    \end{array}
 \right.\,,
\end{equation}
cf. \cite[Thm.2.4]{Ra89} resp.
\cite[p.104]{Hi84}.

\bigskip

{\sc Case 1.1:} Assume in addition that 
$\ind(c^2)\equiv \ind(c)\pmod{2}\,,$ 
i.e.
 $\ind(c^m)\equiv n-1\pmod{2}$ for all $m\ge 1.$
 Then Equations~\eqref{eq:morse}, 
 \eqref{eq:ceins} and \eqref{eq:monotone}
 imply for $k \le N_1:$
 $$w_k=\#\{m \in \n\,;\, \ind(c^m)=k\}=b_k,$$
 and
 \begin{eqnarray*}
 \left(\ind(c^m)\right)_{m\le m_1}=\left(n-1,n+1,\ldots,
 3n-5, 3n-3,3n-3,3n-1, \ldots,\right.
 \\
 \left.
 \ldots, 5n-7,5n-5,5n-5,5n-3,\ldots\right)\,.
 \end{eqnarray*}
In particular we obtain for all $m\le m_1:$
\begin{equation}
\label{eq:indcm}
 \ind(c^{nm})=(2m+1)(n-1)\,.
\end{equation}
Since $\ind(c^{2n})=5n-5 <N_1$ we obtain $m_1\ge 2.$
And Equation~\eqref{eq:indcm} implies
$$
\ind(c^{nm})-m \,\ind(c^n)=
-(m-1)(n-1)$$
for all $m\le m_1,$ hence
Equation~\eqref{eq:minimalequality} 
holds for $r=n.$

\bigskip

{\sc Case 1.2:} Assume in addition that 
$\ind(c^2)\equiv\ind(c)+1\pmod{2},$ i.e.
$\ind(c^m)\equiv m \pmod{2}.$
Then Equations~\eqref{eq:morse}, 
 \eqref{eq:ceins} and \eqref{eq:monotone}
 imply for $k \le N_1:$
$$w_k=\#\{l \in \n\,;\, \ind(c^{2l-1})=k\}=b_k,$$
and
$\ind(c^{2n-1})=\ind(c^{2n+1})=3n-3$ contradicting
Inequality~\eqref{eq:monotone1}.
Therefore this case cannot occur.
 
\bigskip
 
 {\sc Case 2:} Assume $n$ odd, 
then
 \begin{equation}
  b_k=
  \left\{
  \begin{array}{ccl}
  2 &;& k=j (n-1), j\ge 2\\
  1&;& k\ge n-1, k \mbox{ even }, k \not=j (n-1), j\ge 2\\
  0&;& \mbox{ otherwise }
    \end{array}
 \right.\,,
\end{equation}
cf. \cite[Thm.2.4]{Ra89} resp.
\cite[p.104]{Hi84}.

\bigskip

{\sc Case 2.1:} Aussume in addition that 
$\ind(c^2)\equiv \ind(c) \pmod{2},$ 

i.e.
 $\ind(c^m)\equiv 0 \pmod{2}$ for all $m\ge 1.$
  Then Equations~\eqref{eq:morse}, 
  \eqref{eq:ceins} and \eqref{eq:monotone}
  imply 
  for $k \le N_1:$
 $$w_k=\#\{l \in \n\,;\, \ind(c^l)=k\}=b_k,$$
 i.e.
 \begin{eqnarray*}
 \left(\ind(c^m)\right)_{m\le m_1}=\left(n-1,n+1,\ldots,
2n-4, 2n-2, 2n-2, 2n, \ldots,\right.
 \\
 \left.
 \ldots, 3n-5, 3n-3, 3n-3, 3n-1,\ldots\right)\,.
 \end{eqnarray*}
In particular we obtain for all $m\le m_1$
\begin{equation}
\label{eq:indcma}
 \ind(c^{m(n+1)/2})=(m+1)(n-1)\,.
\end{equation}
Since $\ind(c^{2(n+1)})=5n-5$ we obtain $m_1\ge 4.$
Equation~\eqref{eq:indcma} implies
$$
\ind(c^{m(n+1)/2})-m \,\ind(c^{(n+1)/2})=
-(m-1)(n-1)$$
for all $m\le m_1\,,$ 
which is
Equation~\eqref{eq:minimalequality} 
for $r=(n+1)/2.$

\bigskip

{\sc Case 2.2:} Assume in addition that 
$\ind(c^2)\equiv n$ i.e.
$\ind(c^m)\equiv m \pmod{2}.$
Then Equations~\eqref{eq:morse}, 
 \eqref{eq:ceins} and \eqref{eq:monotone}
 imply for $k \le N_1:$
$$w_k=\#\{l \in \n\,;\, \ind(c^{2l-1})=k\}=b_k,$$
and
$\ind(c^{n})=\ind(c^{n+2})=2n-2$ contradicting
Inequality~\eqref{eq:monotone1}.
Therefore this case cannot occur.

  In Case 1.1 and Case 2.1 $m_1\ge 2$ and 
   $\ind(c^{2r})-2\,\ind(c^r)=-(n-1)\,.$ 
   Therefore it follows
  from~\cite[Lemma 3.1(ii)]{BTZ} that 
  $c$ is of elliptic-parabolic type. 
  And the linearized Poincar\'e mapping decomposes into
  $(n-1)$ $(2 \times 2)$ blocks conjugate to rotations.
     \qed
  
\section{Example}
\label{sec:example}
For $\epsilon>0, \eta>0$ with
$\eta+\epsilon<1/2$ let 
$h=h_{\epsilon,\eta}:\R \longrightarrow [0,1]$
be a smooth funtion satisfying
$$
h_{\epsilon, \eta}(t)=\left\{
\begin{array}{ccc}
1 &;& t\in [-\eta,\eta]\\
0 &;& |t|\ge \eta+\epsilon/4
\end{array}
\right.
$$
Let $H_{\epsilon}: \R \longrightarrow \R$
be the smooth function with 
$H_{\epsilon}(t)=\epsilon h_{\epsilon,\epsilon/4}(t).$
And
let $a=a_{\epsilon, \eta}:\R \longrightarrow [-1,1]$
be the smooth function with
$a_{\epsilon,\eta}(t+2)=-a_{\epsilon,\eta}(t)$ and 
$a_{\epsilon,\eta}(1+t)=h_{\epsilon, \eta}(t)$
for all $-1\le t\le 1\,.$

Let $p_1,p_2,\ldots, p_n$ be numbers which are relatively prime
and let $p=p_1p_2\cdots p_n.$
Let $e_0,e_1,e_2,\ldots,e_{2n}$ be an 
oriented 
orthonormal basis of
$\R^{2n+1}.$ 

For the standard Riemannian metric on $S^{2n}\subset \R^{2n+1}$ 
an isometric $\sone$-action 
$\phi: S^{2n} \times \sone  \longrightarrow S^{2n}$
is defined by $\phi\left(e_0,t\right)=e_0$ and
\begin{eqnarray}
\label{eq:phi}
\phi\left(e_{2j-1},t\right)=&\cos(2 \pi p t/p_j)e_{2j-1}
+\sin(2 \pi p t/p_j)e_{2j} \\
\phi\left(e_{2j},t\right)=&-\sin(2 \pi p t/p_j)e_{2j-1}
+\cos(2\pi p t/p_j)e_{2j}\,.\nonumber
\end{eqnarray}
with $j=1,2,\ldots,n.$
We assume that there is a non-reversible Finsler metric
$f$ invariant under $\phi$
for which there are only finitely many geometrically
distinct closed geodesics $c_1,c_1^{-1},$ $c_2,c_2^{-1},\ldots,
c_n,c_n^{-1}$ which occur in pairs differing only
by orientation. 

The geodesic $c_j, j=1,2,\ldots,n$ is the intersection of the
$2$-plane generated by $e_{2j-1}, e_{2j}$ 
invariant under $\phi,$
i.e. $c_j(t)=\cos(2\pi t)e_{2j-1}+
\sin(2\pi t) e_{2j}.$
The closed geodesics $c_j, j=1,2,\ldots,n$ 
do not intersect pairwise.

This assumptions are satisfied by the Katok metrics
$N_{\alpha}$ on $S^{2n}$ with
irrational $\alpha\in(0,1)$ discussed in detail by 
Ziller~\cite[Sec.1]{Zi}. Here the lengths
$L(c_j), L(c_j^{-1})$ of the closed geodesics 
$c_j,c_j^{-1}$ depend on the parameter
$\alpha.$ 

\begin{proposition}
Let $f$ be a bumpy
Katok metric on $S^{2n}$ invariant under
the $\sone$-action $\phi$  by Equation~\eqref{eq:phi}
with only $2n$ geometrically distinct closed geodesics
$c_1,c_1^{-1},\ldots, c_n,c_n^{-1}$
and let 
$L_1:=\max\{L(c_j),L(c_j^{-1}), j=1,2,\ldots,n\}.$
The $n$ closed geodesics $c_1,\ldots,c_n$ 
are invariant under $\phi$ and do not
intersect pairwise.

Fix $k\ge 4.$ For any $L>L_1$ 
there is a non-reversible
Finsler metric $F$ invariant under $\phi$ 
of positive flag curvature and $C^k$-arbitrarily
close to $f$ with only $2n$ geometrically distinct
closed geodesics $d_1,d_2,\ldots,d_{2n}$
of length $<L.$ These $2n$ closed geodesics 
are invariant under $\phi$ and
do not intersect pairwise.
\end{proposition}

Fix the closed geodesic $c=c_1$
and $m=m_1=p/p_1.$ 
Then there is a $\sone$-invariant tubular
neighborhood $U(c)$ of the closed geodesic $c_1$
with coordinates 
$$ (x,t) \in D^{2n-1}_{\epsilon} \times_{\z_{m}} \sone
\longmapsto 
\phi\left(
\cos\left(\|x\|\right)e_1+\sin\left(\|x\|\right)\frac{x}{\|x\|}\,,
\,t\right)\in U(c)
$$
Here $D^{2n-1}_{\epsilon}=\{x\in \R^{2n-1};\|x\|<\epsilon\}$
is a disc of radius $\epsilon$ 
on the Euclidean space generated by $e_0,e_3,e_4,\ldots,e_{2n}$
which we identify via the exponential map of the standard
Riemannian metric on $S^{2n}$
with an open neigborhood of $c(0)=e_1$ on the totally geodesic
hypersphere through $c(0)$ orthogonal to the $\sone$-action.

The $\sone$-action induces
an isometric $\z_m$-action on $D^{2n-1}_{\epsilon},$
and therefore a diagonal action on the
product $D^{2n-1}_{\epsilon}\times \sone$
with quotient space $D^{2n-1}_{\epsilon} \times_{\z_m} \sone.$
Then there are induced coordinates
$(x,t,\xi,\tau)\in
\left(D^{2n-1}_{\epsilon} \times_{\z_m} \sone \right)
\times \R^{2n-1}\times \R$
on the tangent bundle $TU(c)$ restricted to the tubular
neighborhood $U(c)$ of $c.$
Because of the $\sone$-symmetry the Finsler metric $f$
in the coordinates $(x,t,\xi,\tau)$ does not depend 
on $t,$ i.e.
$f\left(x,t,\xi,\tau\right)=f\left(x,\xi,\tau\right).$

Choose $\eta \in (0,1/4)$ and a sufficiently small 
$\epsilon>0.$ For a sufficiently small $s>0$ the mapping
$\Psi=\Psi_{s,\epsilon,\eta}:
TU(c)-T^0U(c)\longrightarrow TU(c)-T^0U(c)$
of the tangent bundle $TU(c)$ 
of the tubular neighborhood
minus the zero section
$T^0U(c)$ defined by
$$\Psi\left(x,t,\xi,\tau\right)=
\left(x+s\, H_{\epsilon}\left(\|x\|\right)\,
a_{\eta,\epsilon}\left(\frac{\tau}{\sqrt{\tau^2+\|\xi\|^2}}\right)
\,e_0,
t,\xi,\tau\right)$$
is a diffeomorphism.
For a fixed $k \ge 4$ one can choose $s>0$ sufficiently small
such that 
the diffeomorphism $\Psi_{s,\epsilon,\eta}$ is
arbitrarily close to the identity map in
the $C^k$-norm.
The mapping is well-defined since $e_0$ is fixed under
$\phi.$

The diffeomorphism 
$\Psi$
extends by the identity onto
the complement of $TU(c)$ on $TS^{2n}.$
And the diffeomorphism $\Psi$ is positively homogeneous
of degree $1,$
i.e.
$\Psi\left(x,t,\lambda \xi, \lambda \tau\right)
=\lambda \Psi\left(x,t,\xi, \tau\right)$
for $\lambda > 0.$

Then define the Finsler metric
$$\overline{f}(x,t,\xi,\tau)=
\overline{f}_{s,\epsilon,\eta}(x,t,\xi,\tau)=
f\left(\Psi^{-1}_{s,\epsilon,\eta}(x,t,\xi,\tau)\right)\,.$$
Let
$V_{\pm}=\left\{(x,t,\xi,\tau); \|x\|<\epsilon/4,
\pm \tau/\sqrt{\tau^2+\|\xi\|^2}> 1-\eta
\right\}$
and \\
$V_1=\left\{(x,t,\xi,\tau); 
|\tau|/\sqrt{\tau^2+\|\xi\|^2} \le 1-\eta-\epsilon/4
\right\}.$
Then
\begin{equation}
\label{eq:vplusminus}
\overline{f}_{s,\epsilon, \eta}(x,t,\xi,\tau)=
\left\{
\begin{array}{ccc}
f(x,t,\xi,\tau) &;& (x,t,\xi,\tau) \in V_1\\
f\left(x\mp s \epsilon e_0,t,\xi,\tau\right)
&;& (x,t,\xi,\tau) \in V_{\pm}
\end{array}
\right. \,.
\end{equation}
Therefore we obtain
that $d_1: t \mapsto (s \epsilon e_0,t)$ 
and $d_2: t \mapsto (-s \epsilon e_0,-t)$
are two closed geodesics
of $\overline{f}_{s,\epsilon,\eta}$ 
in $U(c)$ which do not intersect.

The flag curvature of the Katok metric is constant with
value $1,$ 
cf.~\cite[Sec.4]{Ra04}. 
And we obtain that the flag curvature 
$\overline{K}\left(y,Y,\sigma\right)$
of a flag 
$(y,Y,\sigma)$ with 
$y \in S^{2n}, Y \in T_yS^{2n}$ and $\sigma$
a two-dimensional subspace of $T_yS^{2n}$
 containing
$y$ with respect to the perturbed Finsler metric
$\overline{f}=\overline{f}_{s,\eta,\epsilon}$
can be different from $1$ only if
$y \in U(c)$ and if
the coordinates $(x,t,\xi,\tau)$ of the tangent vector
$Y$ satisfy: $(x,t,\xi,\tau) \not\in V_1\cup V_+ \cup V_-.$
The Finsler metric
$\overline{f}_{s,\epsilon,\eta}$ is invariant under $\phi.$
But the symmetry with respect to the reflection at
the $(e_1,e_2)$-plane is broken.

We can use this perturbation 
inductively in tubular neighborhoods 
for all closed geodesics $c_2,\ldots,c_n$
and obtain 
a non-reversible Finsler metric 
$F=F_{s,\epsilon,\eta}$
on $S^{2n}$
invariant under $\phi$ 
with closed geodesics
$d_3, d_4,\ldots,d_{2n}$ 
also invariant under $\phi.$ Any two of these closed
geodesics do not intersect.

Now fix $L >L_1\,.$ Since 
$f$ is bumpy there is a sufficiently small
$s>0$ such that all
closed geodesics 
of the Finsler metric $F_{s,\epsilon,\eta}$
of length $<L$ are geometrically
equivalent to $d_1, d_2,\ldots, d_{2n},$
for this argument compare for example
~\cite[p.12, 13]{An}.
Here we use the strong $C^k$-topology with
$k\ge 4$ (instead of $k \ge 2$ in
the Riemannian case) since in contrast to the Riemannian
case the geodesic coefficients depend on
fourth derivatives of the Finsler metric.

Since 
by Equation~\eqref{eq:vplusminus}
the Finsler metrics $f$ and $F=F_{s,\epsilon,\eta}$
are isometric in open neighbhorhoods of the 
velocity fields $c_j'$ and $d_{2j-1}'$ resp.
$-c_j'$ and $d_{2j}'$ 
for $j=1,2,\ldots,n$
the flag curvature is constant in these
neighborhoods and the lengths 
$L(c_j)=L(d_{2j-1}), L(c_j^{-1})=L(d_{2j}), j=1,2,\ldots,n$
coincide. Therefore the 
picture from Morse theory of $f$ 
produced
by $c_1,c_1^{-1},\ldots,c_n,c_n^{-1}$ 
and $F$ produced by
$d_1,d_2,\ldots,c_{2n}$ up to length $L$ coincide.

For $n=2$ it follows from~\cite[Sec.6]{HP}
that also the perturbed metric $F$ on $S^2$
is dynamically convex.
Therefore there are either infinitely many geometrically
distinct closed geodesics or there are only two geometrically
distinct closed geodesics.

For $n=2$ we obtain non-reversible and rotationally
invariant Finsler metrics
in any neighborhood of the Katok metric with two
closed geodesics $d_1,d_2$ 
also invariant under the $\sone$-action which do not intersect.
For existence results for closed geodesics invariant under
isometries cf.~\cite[Prop.2, Prop.3]{Ra07}.


\end{document}